
\input amstex
\define\Date{July 25, 2001}
\documentstyle{amsppt}
\magnification=1200
\loadeurm
\loadeusm
\loadbold
\font\twbf=cmbx10 scaled\magstep2
\font\chbf=cmbx10 scaled\magstep1
\font\hdbf=cmbx10


\pagewidth{6.5 true in}
\pageheight{8.9 true in}
\nologo
\overfullrule=3pt
\define\0{^{@,@,\circ}}
\define\1{^{-1}}
\define\2{^{\hbox{\rm h}}}
\define\8{^{\infty}}

\define\ansp{\eurm a}

\define\blank{\underline{\phantom{J}}}

\define\cJ{{\Cal J}}

\define\clo#1#2#3{{{#1}^{#2}}_{#3}}

\define\cS{{\Cal S}}
\define\cT{{\Cal T}}
\define\cln{\colon}
\define\col#1{\,\colon_{#1}}
\define\Col{\,\colon\,}

\define\demb#1{\demo{\bf #1}}
\define\depth{\text{depth}\,}

\define\er{^{(e)}R}
\define\F{\boldkey F}
\define\gr{\hbox{\rm gr}}

\define\Ic{\widehat{I}}

\define\inc{\subseteq}
\define\inj{\hookrightarrow}
\define\intp #1 {\lfloor #1 \rfloor}

\define\m{\bold m}

\define\N{{\Bbb N}}

\def\pa{\parindent 35 pt}
\define\part#1{\item{\hbox{\rm (#1)}}}
\define\pd{\hbox{\rm pd}}
\define\pl{^\boldkey{+}}

\define\Pnl{\boldkey{(}}
\define\Pnr{\boldkey{)}}
\define\q{^{[q]}}
\define\qq{^{1/q}}

\define\Rc{\widehat{R}}
\def\sy#1 {^{\Pnl#1\Pnr}}

\define\Spec{\text{Spec}\,}

\define\Tor{\hbox{\rm Tor}}
\define\ux{\underline{x}}
\define\uy{\underline{y}}
\define\vect#1#2{#1_1,\, \ldots, #1_{#2}}
\define\Z{{\Bbb Z}}
\quad\bigskip\bigskip
\centerline{\twbf COMPARISON OF}
\vskip 10 pt plus .5 pt minus .5 pt
\centerline{\twbf SYMBOLIC AND ORDINARY}
\vskip 10 pt plus .5 pt minus .5 pt
\centerline{\twbf POWERS OF IDEALS}

\topmatter
\author
by Melvin Hochster and Craig Huneke
\endauthor
\rightheadtext{COMPARISON OF SYMBOLIC AND ORDINARY POWERS OF IDEALS}
\leftheadtext{Melvin Hochster and Craig Huneke}
\thanks{The authors were supported in part by grants from the National
Science Foundation.} \endthanks
\thanks{Version of \Date.}\endthanks
\endtopmatter

\document
{\baselineskip = 15.95 pt
\parindent = 12 pt
\quad\bigskip
{\chbf
\centerline{\hdbf 1. Introduction}
}
\bigskip

All given rings in this paper are commutative, associative with
identity, and Noetherian.
Recently,  L.~Ein, R.~Lazarsfeld, and K.~Smith [ELS] discovered a
remarkable and surprising fact about the behavior of symbolic
powers of ideals in affine regular rings of equal characteristic 0:
if  $h$ is the largest height of an associated prime of $I$,
then $I\sy hn \inc I^n$ for all $n \geq 0$.  Here,
if $W$ is the complement of the union of the associated primes of
$I$,  $I\sy t $ denotes the contraction of $I^t R_W$ to $R$,  where
$R_W$ is the localization of $R$ at the multiplicative system $W$.
Their proof depends on the theory of multiplier ideals, including an
asymptotic version, and, in particular, requires resolution of singularities
as well as vanishing theorems.
We want to acknowledge that without their generosity and quickness in
sharing their research this manuscript would not exist.

Our objective here is to give stronger results that can be proved by
methods that are, in some ways, more elementary.
Our results are valid in both equal
characteristic 0 and in positive prime characteristic $p$, but
depend on reduction to characteristic $p$.
We use tight closure methods and, in consequence, we  need neither
resolution of singularities nor vanishing
theorems that may fail in positive characteristic.  For the most
basic form of the result, all that we need from tight closure theory
is the definition of tight closure and the fact that in a regular ring,
every ideal is tightly closed.  We note that the main argument here is
closely related to a proof given in [Hu, 5.14--16, p.\ 45] that regular
local rings in characteristic $p$ are UFDs, which proceeds
by showing that Frobenius powers of height one primes
are symbolic powers.

Our main results in all characteristics are summarized in the following
theorem. Note that $I^*$ denotes the tight closure of
the ideal $I$. The characteristic zero notion of tight closure used in
this paper is the {\it equational tight closure} of [HH6] (see, in
particular Definition (3.4.3) and the remarks in (3.4.4) of
[HH6]).  This is the smallest
of the characteristic zero notions of tight closure, and therefore
gives the strongest result.   See \S3.1 for
a discussion of the Jacobian ideal $\cJ(R/K)$ utilized in part (c).

\proclaim{Theorem 1.1} Let $R$ be a Noetherian ring containing
a field.  Let $I$ be any  ideal of
$R$, and let $h$ be the largest height\footnote{The results stated here
are all valid if one defines $h$ instead to be the largest
analytic spread of $IR_P$ for any associated prime $P$ of $I$, which, in
general, may be smaller:  see Discussion 2.3} of any associated
prime of $I$.
{\pa
\part{a} If $R$ is regular, $I\sy hn+kn \inc (I\sy k+1 )^n$
for all positive $n$ and  nonnegative $k$.  In particular,
 $I\sy hn \inc I^n$ for all positive integers $n$.

\part{b}
\parindent = 12pt
If $I$ has finite projective
 dimension
  then $I\sy hn \inc (I^n)^*$ for all positive integers $n$. \par
\pa
\part{c} If $R$ is finitely generated, geometrically reduced
(in characteristic 0,
this simply means that $R$ is reduced) and equidimensional over
a field $K$, and locally $I$ is either $0$ or contains a nonzerodivisor
(this is automatic if $R$ is a domain),
then, with $J = \cJ(R/K)$,  for every nonnegative integer $k$ and
positive integer $n$, we have that $J^nI\sy hn+kn \inc ((I\sy k+1 )^n)^*$ and
$J^{n+1}I\sy hn+kn \inc (I\sy k+1 )^n$.   In particular, we have that
$J^nI\sy hn \inc (I^n)^*$ and $J^{n+1}I\sy hn  \inc I^n$ for all positive
integers $n$.
\par
 }
\endproclaim

These results, when specialized to the case where $R$ is
regular, recover the cited result from [ELS].

The theorem above is a composite of Theorems 2.6, 3.7, and 4.4
below.

We note that by results\footnote{E.g., it is shown in [Swsn] that if
$I\inc J$ are ideals of a Noetherian ring, and we let
$I\Col J\8 = \bigcup_t I\Col J^t$,
then if the
$I$-adic filtration is equivalent to the $I^n \Col J\8$ filtration,
there exists an integer $h'$ such that for all $n$,
$I^{h'n} \Col J\8 \inc I^n$.} of [Swsn] one expects, in many cases, to have
results that assert that, given a fixed ideal $I$ in
a Noetherian ring,  for some choice of positive integer $h'$
(independent
of $n$ but depending on $I$) one has
$I\sy h'n \inc I^n$  for all positive integers $n$.
 What is not expected is
the very simple choice of $h'$ that one can make in a regular ring, and
the extent to which it is independent of information about $I$.  E.g.,
when $d = \dim R$ is finite, then with $a = d$
one has that $I\sy an \inc I^n$ for all
ideals $I$ (if $R$ is local  one has this for all unmixed $I$ with
$a = \hbox{max}\,\{d-1,\,1\}$ ---
one does not have to worry about letting $h = d$,
since for the maximal ideal one has that ordinary and symbolic powers
coincide).

We conclude this introduction by sketching the proof
of Theorem 1.1(a) for regular domains in characteristic $p > 0$
when $k=0$  in the special case where $I$ is a radical ideal.
The proof is very simple and brief
in this case, and we hope that this argument will help
the reader through the complexities of the rest of this paper.
Suppose that $I \not= (0)$ is a radical ideal, and let $h$ be the largest
height of any minimal prime.  If $u\in I\sy hn $, then for
every $q = p^e$ we can write $q = an+r$ where $a \geq 0$ and
$0 \leq r \leq n-1$ are integers.  Then $u^a \in I\sy han $
and $I^{hn} u^a \inc I^{hr}u^a \inc I\sy han+hr = I\sy hq $.
We now come to a key point:  we can show that ($*$) $I\sy hq \inc I\q$.
To see this, note that because the Frobenius endomorphism
is flat for regular rings, $I\q$ has no associated primes other
than the minimal primes of $I$ (cf.\ Lemma (2.2d)), and it suffices
to check ($*$)  after localizing at each minimal prime $P$ of $I$.
But after localization, $I$ has at most $h$ generators, and so each
monomial of degree $hq$ in these generators is a multiple
of the $q\,$th power of at least one of the generators.
This completes the proof of ($*$).
Taking $n\,$th powers gives that $I^{hn^2}u^{an} \inc (I\q)^n =
(I^n)\q$, and since $q \geq an$, we have that $I^{hn^2}u^q \inc (I^n)\q$
for fixed $h$ and $n$ and all $q$.  Let $d$ be any nonzero element
of $I^{hn^2}$.  The condition that $du^q \in (I^n)\q$ for all $q$
says precisely that $u$ is in the tight closure of $I^n$ in $R$.
But in a regular ring, every ideal is tightly closed
(cf.\ [HH2, Th. (4.4)]), and so $u \in I^n$,  as required. \qed
\bigskip
\bigskip
{\chbf
\centerline{\hdbf 2. The regular case in characteristic $\boldkey{p}$}
}
\bigskip

\demb{Discussion 2.1} We recall
some terminology and notation.  $R\0$ denotes the complement of
the union of the minimal primes of $R$, and so, if $R$ is reduced,
$R\0$ is simply the multiplicative system of all nonzerodivisors in $R$.
We shall write $\F^e$ (or $\F^e_R$ if we need to specify the base ring)
for the {\it Peskine-Szpiro} or {\it Frobenius} functor from $R$-modules
to $R$-modules.  This is a special case of the base change functor
from $R$-modules to $S$-modules that is simply given by
$S \otimes_R\,\blank$:
in the case of $\F^e$, the ring $S$ is $R$, but the map $R \to R$
that is used for the algebra structure is the $e\,$th iteration $F^e$
of the Frobenius endomorphism:  $F^e(r) = r^{p^e}$.
We shall use the notation $\er$ for $R$ viewed as an $R$-algebra
via the homomorphism $F_R^e\cln R \to R$.

In particular, if $M$ is given as the
cokernel of the map represented by a matrix $\bigl(r_{ij}\bigr)$,
then $\F^e(M)$ is the cokernel of the map represented by the
matrix $\bigl(r_{ij}^{p^e}\bigr)$.  Unless otherwise indicated, $q$
denotes $p^e$ where $e \in \N$.  For $q = p^e$,
$\F^e(R/I) \cong R/I\q$,  where $I\q$ denotes the ideal generated
by the $q\,$th powers of all elements of $I$ (equivalently, of
generators  of $I$).  Note that $\F^e$ preserves both freeness and
finite generation of modules, and is exact precisely when
$R$ is regular (cf. [Her], [Kunz]).  \enddemo

\proclaim{Lemma 2.2 (Peskine-Szpiro)} Let $R$ be a
Noetherian ring of characteristic $p$,
and $M$ be a finitely generated $R$-module of finite projective
dimension over $R$.
Then:
{\pa
\part{a} For all $i \geq 1$,  $\Tor^R_i(M,\,\er) = 0$.
\part{b} If $R$ is local and one applies $\F^e$ to a minimal
free resolution of $M$,  one obtains a minimal free
resolution of $\F^e(M)$. In particular, $\pd_RM = \pd_R\F^e(M)$.
\part{c} For all $e \geq 1$,  the set of associated primes of
$M$ is the same as the set of associated primes of $\F^e(M)$.
\part{d} In particular, if $R$ is regular, so that $\F^e$ is flat,
then the conclusions of {\rm (b)} and {\rm (c)} are
valid for every finitely generated
$R$-module. \par }
\endproclaim
\demb{Proof} We refer to [PS] for part (a).
Part (b) is well-known and is immediate from (a). Part (c) is likewise
well-known,
but we mention that it reduces to the local case by
localization at a given prime, and so it reduces to checking that
the maximal ideal of $R$ is associated to $M$ if and only if
it is associated to $\F^e(M)$.  But when $M$ is a module of
finite projective dimension,  $m$ is associated to $M$
if and only if $\pd_R M = \depth R$.

We refer the reader to [PS], [Her], and [Kunz] for related results.
\qed \enddemo

\demb{Discussion 2.3: integral dependence of ideals, analytic spread,
and minimal reductions}  Recall that an element $r$  of a ring
$R$ is {\it integrally dependent} on an ideal $I$ if there is an
integer $t \geq 1$ and an equation of the form
$r^t + i_1r^{t-1} +\,\cdots\,+i_kr^k+\,\cdots\,+i_{t-1}r + i_t = 0$,
where  $i_k \in I^k$,  $1 \leq k \leq t$.  The elements of the
ring $R$ integrally
dependent on $I$ form an ideal $J \supseteq I$,
the {\it integral closure} of
$I$ in $R$.  We refer the reader to [L1],
\S5 of [HH2], and [NR1--2] for more detailed information about integral
dependence and analytic spread.

(a) In a Noetherian local ring $(R,\,m,\,K)$ with maximal ideal
$m$ and residue field $K$,  the {\it analytic spread} $\ansp (I)$ of
an ideal
$I \inc m$ is the Krull dimension of the ring
$$
K \otimes_R \gr_I R \cong
K \oplus I/mI \oplus I^2/mI^2 \oplus\, \cdots \oplus I^k/mI^k \oplus \cdots ,
$$
The analytic spread is a lower bound on the least number of generators
of an ideal $I_0 \inc I$ such that $I$ is integrally dependent on $I_0$.
If $K$ is infinite, there always is an ideal $I_0$ with $\ansp(I)$
generators such that $I$ is integral over $I_0$, and such an ideal
$I_0$ is called a {\it minimal reduction} of $I$.  When $K$ is
infinite, one may find
generators for a minimal reduction  $I_0$ of $I$ by simply
taking a linear homogeneous system of parameters for $K \otimes_R \gr_I R$,
say $\vect f \ansp \in I/mI \cong [K \otimes_R \gr_I R]_1$,
and lifting the  $f_j$ to elements of $I$.

Note that the analytic spread of $I$ is bounded both by the number
of generators of $I$ and by the Krull dimension of $R$.

(b) If $I$ is an ideal of $R$ and $t$ is an indeterminate over $R$,
then the associated primes of $IR[t]$ are those of the form
$Q = PR[t]$ where $P$ is an associated prime of $I$.  For this
$Q$ the analytic spread of $IR[t]_Q$ is the same as the analytic
spread of $IR_P$.  Thus, the maximum analytic spread after localization
at an associated prime is the same for $IR[t]$ in $R[t]$ as it is
for $I$ in $R$.    Moreover, the symbolic powers of $IR[t]$ are
the expansions of the symbolic powers of $I$.

(c) If  $S$ is flat over $R$ then the maximum analytic spread of $IS$
after localizing at an associated prime in  $S$ is at most what
it was for $I$ in $R$.  To see this, first note that by replacing
$R \to S$
by $R[t] \to S[t]$ and $I$ by $IR[t]$, we may assume without
loss of generality that
the residue fields of the local rings of associated primes of $I$
in $R$ are infinite.  We return to the original notation.
By Proposition 15 in Section IV B.4. of
[Se],  $Q$ is an associated prime of $IS$ if and only if it is
an associated prime of (0) in $S/PS$ for some associated prime
$P$ of $I$.    But $S/PS$ is flat and, hence, torsion-free over
the domain $R/P$, which implies that $Q$ lies over $P$.  Thus,
we have a map  $R_P \to S_Q$.  If the analytic spread of
$IR_P$ is $h$, it is integral over an ideal with $h$ generators.
But then $IS_Q$ is integral over the expansion of the same ideal,
and the result follows.

(d) We recall also that in a Noetherian ring $R$,  $I$ is integrally
dependent  on $I_0$ if and only if there exists an integer $k$ such for all
positive integers $n$, $I^{k+n} =
I^kI_0^n$.  In particular, it then follows that $I^{k+n} \inc I_0^n$
for all positive integers $n$.
\enddemo
\medskip
Part (b) of the next result plays a critical role in the proofs
of our theorems.  It is closely related to the Brian\c con-Skoda
theorem, and related results were used to prove a tight closure
form of the Brian\c con-Skoda theorem in \S5 of [HH2].  In fact,
our first proofs of some of the results here made use of the
Brian\c con-Skoda theorem in a sharpened form\footnote{To be precise,
if $I$ is integral over an ideal with at most
$h$ generators in a Noetherian domain of
characteristic $p > 0$, then the integral closure of
$I^h$ is contained in $I^+$. Here,  if $R\pl$ denotes
the integral closure of $R$ is an algebraic closure of its
fraction field (which is unique up to non-unique
isomorphism), then $I^+ = IR\pl \cap R$. It is known that $I \inc I^+ \inc
I^*$, while equality for ideals generated by part of a system of
parameters in an excellent local domain is established in [Sm1].
Because the formation of $R\pl$ commutes with localization at
a multiplicative system, plus closure commutes with localization.
We refer the reader to [HH4--5] and [Sm1] for further discussion
of $R\pl$. }
given in [HH5], Theorem (7.1), that uses plus closure instead
of tight closure, together with the fact that
plus closure commutes with localization.

\proclaim{Lemma 2.4}
 Let $R$ be a ring.
{\pa
\part{a}
If $I = (\vect u h)$ then for all integers $t \geq 1$ and $k \geq 0$,
 $$I^{ht+kt -h+1} \inc (u_1^t,\, \ldots,\, u_h^t)^{k+1}.$$
In particular,
$I^{ht+kt} \inc (u_1^t,\, \ldots,\, u_h^t)^{k+1}$. Hence,
if $R$ has prime characteristic $p> 0$ and $q = p^e$ is a power
of $p$,  then  $I^{hq+kq} \inc (I\q)^{k+1} = (I^{k+1})\q$.

\part{b} (Key Lemma) Let $R$ be Noetherian of positive
prime characteristic $p$.   Suppose that $I$ is an ideal
of $R$, that $W$ is the complement of the union of the
associated primes of $I$,  and that ${}^W$ indicates the result
of expanding an ideal of $R$ to $R_W$ and then contracting it to $R$.
Suppose that for every associated prime $P$ of $I$,
$IR_P$ has analytic spread at most  $h$ in  $R_P$.  Then
there is a fixed positive integer $s$ (depending on $I$) with the following
property:\smallskip
\item{}For all choices of integers $n \geq 0$,  $q = p^e$, and
$k \geq 0$, we have that if
$u \in I\sy hn+kn $ then
$$
I^{s + (h+k)(n-1)}u^{\intp q/n } \in \bigl((I\sy k+1 )\q\bigr)^W ,
$$
where $\intp q/n $ denotes the integer part of $q/n$.
If $R$ is regular or if $I$ has finite projective dimension
and $k=0$,  the superscript $^W$ can be omitted.  \par
}
\endproclaim
\demb{Proof}  Consider any monomial
$u_1^{b_1}\,\cdots \,u_h^{b_h}$ in the
$u_j$ in which the sum of the exponents,  $b_1 +\,\cdots\,+b_h$,
is at least $ht+kt -h+1$.  Write each $b_j = a_jt + c_j$ where
$a_j$ is a positive integer and $0 \leq c_j \leq t-1$.  Then it
suffices to prove that the sum of the $a_j$ at least $k+1$, for then
the original monomial is a multiple of $(u_1^t)^{a_1}\,\cdots\,(u_h^t)^{a_h}
\in (u_1^t,\,\ldots,\,u_h^t)^{k+1}$.  But, otherwise, we have that
the sum of the $a_k$ is at most $k$,  which means that the sum of
the $b_j$ is at most $kt+h(t-1) < ht+kt-h+1$, a contradiction.
(A slightly weaker version is proved in a parenthetical comment
near the bottom of p.~45 of [HH2].) The remaining statements in
(a) are immediate.

For part (b), first note that the issues are unaffected by
adjoining an indeterminate $t$ to the ring $R$ and replacing
$I$ by $IR[t]$.  A choice of $s$ that works for $IR[t]$ and $R[t]$
will work for $I$ and the original ring $R$:  the associated primes
of $IR[t]$ are simply those of the form $PR[t]$,  where $P$ is
an associated prime of $I$, and if $W'$ is the complement of
the union of the associated primes of $IR[t]$ we have
that $R[t]_{W'}$ is faithfully flat over  $R_W$.  Moreover,
for every $P$,  the analytic spread of $IR[t]_{PR[t]}$ is the
same as the analytic spread of  $IR_P$.

Thus, we may assume without loss of generality that the residue
field of each of the rings $R_P$ is infinite when $P$ is an associated
prime of $I$, and it follows
that for each associated prime $P_i$ of $I$ we can choose
an ideal $J_i \inc I$ with at most $h$ generators such that
$I_{P_i}$ is integral over $(J_i)_{P_i}$.  By 2.3(d)
there is a positive  integer  $s_i$ such that  $I_{P_i}^{s_i+N}
\inc (J_i)_{P_i}^N$ for all positive integers $N$.  Let  $s$ be the
maximum of the $s_i$.   Write
$q = an +r$ with $0 \leq r \leq n-1$.
To prove that
$I^{s + (h+k)(n-1)}u^a \in
\bigl((I\sy k+1 )\q\bigr)^W$ it suffices to prove that whenever
$P = P_i$ is an associated prime of $I$, we have that
$$
I^{s + (h+k)(n-1)}u^aR_P \inc \bigl((I\sy k+1 )\q\bigr)R_P .
$$
 Let $J = J_i$.
Since $R_P$ contains $W^{-1}$,  after localization at $P$ the
symbolic and ordinary powers of $I$ are the same.  But then (recall that
$q = an +r$ with $0 \leq r \leq n-1$), we have that
$$
I_P^{s+(h+k)(n-1)}u^a\inc  I_P^{ s+(h+k)r +(h+k)an} \inc
I_P^{s + (h+k)(an+r)}\inc J_P^{(h+k)q}
$$
(the last inclusion holds by the choice of $s$).  But
$J_P^{(h+k)q}\inc
(J_P^{k+1})\q$ by part (a), and
this is clearly contained in  $(I_P^{(k+1)})^{[q]})$,  as required.

The omission of $W$ when $R$ is regular is justified by
the fact that $I\sy k+1 $ has the property that no element
of $W$ is a zerodivisor on $R/ I\sy k+1 $, and, since the
Frobenius endomorphism
is flat, the elements of $W$ are also nonzerodivisors
on $R/(I\sy k+1 )\q$ for all $q$ by Lemma 2.2(d). If instead
$k=0$ and $I$ has finite projective dimension we may
apply Lemma 2.2(c) instead.
\qed\enddemo

We shall say that an ideal $I$ of a Noetherian ring $R$
is {\it locally generically free}  if for
every prime ideal $P$ of $R$,  $IR_P$ is
either (0) or contains a nonzerodivisor.  When $R$ is local with
total quotient ring $T$,  this is equivalent to requiring
that $IT$ be free (of rank 0 or 1, necessarily).  Notice that
ideals of finite projective dimension are automatically locally
generically free.  This is well-known, but we indicate a brief
argument.  The point is that in the local case one has a free resolution,
and so whenever one has a localization such that the ideal is
projective ($\equiv$ free), the rank is the same as the alternating
sum of the ranks of the free modules in the finite free resolution.
Once one tensors with the total quotient ring (of the local ring)
one has a semilocal ring with all maximal ideals of depth 0.
By the Auslander-Buchsbaum theorem, all modules of finite projective
dimension are now locally free.  By the remarks above the
rank is constantly 0 or 1,  so that the ideal has become either
$(0)$ or free of rank one.

\proclaim{Lemma 2.5} Let $R$ be a Noetherian ring.
{\pa
\part{a} If $R$ is Noetherian and
$I$ is an ideal containing a nonzerodivisor, then
$I$ is generated by the nonzerodivisors in $I$.
\part{b} If $R$ is Noetherian local and $I \not=0$ has finite projective
dimension, then
it contains a nonzerodivisor.  I.e., if $R$ is any Noetherian ring
and $I$ has finite projective dimension, then $I$ is locally generically
free.
\part{c}
If $R$ is Noetherian with $\Spec R$  connected and
$I \not=0$ is locally generically free, then it contains a nonzerodivisor.
\par
}
\endproclaim

\demb{Proof}  For part (a), let $I_0 \inc I$ be the ideal generated
by all nonzerodivisors in $I$.  Then $I$ is contained in the union
of $I_0$ and the associated primes of $(0)$ in the ring.  Since
$I$ is not contained in any associated prime of (0), we must
have that $I \inc I_0$, and so $I = I_0$.  Part (b) was established in
the discussion preceding the statement
of the lemma.  Finally, to prove (c),  let
$S$ be the set of primes $P$ such
that $IR_P$ contains a nonzerodivisor  and let $T$ be the set of primes
$P$ such that  $IR_P$ is zero. Then $\Spec R$ is the disjoint
union of these two sets.  Both have the property that if $P \inc Q$ and
$Q$ is in the set, then $P$ is in the set.  It follows that if $P \inc Q$
and $P$ is in one of these sets, then $Q$ is in the same set.  Thus,
both sets are Zariski closed.  Since $\Spec R$ is connected, one of these
sets is empty, and since $I \not= (0)$, we have that $T$ is empty.
Then $I$ is not contained in any associated prime $P$ of $(0)$
(or its localization $IR_P$ would consist entirely of zerodivisors).
Hence, there is an element
of $I$ not in any associated prime of $(0)$.
\qed\enddemo

We are now ready to prove one of our main results.

\proclaim{Theorem 2.6} Let $I$ be  ideal of a Noetherian
ring of positive prime characteristic $p$.
Let $h$ be the largest height of any associated
prime of $I$ (or let $h$ be the largest analytic spread of
$IR_P$ for $P$ an associated prime of  $I$).
Then, if $R$ is regular,
$I\sy hn \inc I^n$ for every positive integer
$n$,   while if $I$ has finite projective dimension,
$I\sy hn  \inc (I^n)^*$ for every positive
integer $n$.

If $R$ is regular one has more generally that for every
nonnegative integer $k$,
$I\sy hn+kn \inc (I\sy k+1 )^n$ for every positive integer
$n$. \endproclaim

\demb{Proof}
Since $I$ has finite projective dimension (this is automatic
if the ring is regular),
we may apply Lemma 2.4.  If $R$ is a product we may consider
the problem for the various factors separately, and so we may
assume without loss of generality that $\Spec R$ is connected.
If $I = (0)$ there is nothing to
prove.  Otherwise, by Lemma 2.5(c), $I$ contains a nonzerodivisor.

We handle all cases of the theorem at once by assuming either
that $R$ is regular or that $I$ has finite projective dimension
and that $k=0$.  Choose $s$ as in the Key Lemma 2.4(b).
For every
$q = p^e$ we may write $q = an + r$,  where $a$ is a nonnegative
integer and $0 \leq r \leq n-1$.
Now, $u \in I\sy hn+hk $ implies, by Lemma
2.4(b),  that $I^{s + (h+k)(n-1)}u^a \inc (I\sy k+1 )\q$
(note that in both cases the superscript $^W$ is not needed) and
we may raise both sides to the $n\,$th power to get
$$
I^{sn+(h+k)(n-1)n}u^{an} \inc ((I\sy k+1 )^n)\q .
$$
 We may multiply
by $u^r$, and so abbreviating $b = sn+(h+k)(n-1)n$  we have that
$I^{b}u^q \inc ((I\sy k+1 )^n)\q$
for all $q$.  Since $I$ contains a nonzerodivisor, so does
$I^{b}$: call it $d$.   Notice that $b$, and, hence,  $I^b$, does not
depend on $q$.  We therefore have that $du^q \in ((I\sy k+1 )^n)\q$
for all $q$.  Thus,  $u \in ((I\sy k+1 )^n)^*$.  Since every
ideal is tightly closed in case the ring is regular, the
proof is complete in all cases.
\qed\enddemo

\demb{Remark 2.7}  The result above is also valid for ideals $I$ in
Noetherian rings $R$ if
$V(I)$ is disjoint from the singular locus of $R$, and the
singular locus is closed.  (If  $I_1$ defines the singular
locus then $I + I_1 = R$.  Choose $f \in I_1$ so that it
is a unit in $R/I$.  Then it is also a unit modulo any
any ideal containing a power of $I$.
It follows that $R/I' \cong R_f/I'$ whenever $I' = I^n,\, I\sy hn+kn $, or
$(I\sy k+1 )^n$,  and the result is immediate from this observation
and the fact that we may apply Theorem 2.6 to the regular ring $R_f$.)
Precisely the same observation holds in the equal characteristic
0 case, making use of Theorem 4.4 instead.
 \bigskip
 \bigskip
{\chbf
\centerline{\hdbf
3. Singular affine algebras in positive characteristic} \smallskip
}\bigskip
The results of this section depend heavily on the fact that, in
positive characteristic, the elements of the
Jacobian ideal
can be used not only as test elements, but also have the
property that their $q\,$th powers  ``multiply
away" the effects of embedded components of
$q\,$th bracket powers of unmixed ideals (these
do not occur in the regular case):  a precise statement is
given in Lemma 3.6.  Statements of this sort depend heavily in turn
on the Lipman-Sathaye Jacobian theorem.

\demb{Discussion 3.1:  the Jacobian ideal} Let $A$ be a reduced
Noetherian ring with total quotient ring $T = \cT(A)$, so
that $T$ is a finite product of fields.  Let $R$ be a finitely
generated $A$-algebra such that $R$ is torsion-free
over $A$ (i.e., nonzerodivisors of $A$ are nonzerodivisors on $R$)
and such that $T \otimes_A R$ is equidimensional of dimension
$d$ and geometrically
reduced.   We then
define the {\it Jacobian ideal} $\cJ(R/A)$ as follows.  Choose
a finite presentation of $R$ over $A$,  say
$R \cong A[\vect x n]/(\vect f m)$, and let $\cJ(R/A)$ denote
the ideal of $R$ generated by the images of the
$n-d$ size minors of Jacobian matrix $\bigl(\partial f_j/\partial x_i\bigr)$.
An important case is where $A = K$ is a field.

We note the following easy facts:

{\pa
\item{(1)} $\cJ(R/A)$ is independent of the choice of presentation.
E.g., if one changes the set of generators $f_j$ of the denominator
ideal, it suffices to compare the result from each set of generators
with the union.  By induction, one only needs to see what happens
with one additional generator.  The calculation is then very easy.
If one has two different presentations one can put them together
(think of the two sets of variables as disjoint and independent).
Thus, one need only compare Jacobians when one uses some extra
generators to give a presentation, and, by induction, it suffices
to consider the case of one extra generator.  But then the
denominator ideal has the form $\vect f m, y-g(\vect x n)$ where
$y$ is a new variable and
$g$ maps to the extra generator in $R$.  Again, the calculation
is now easy.

\item{(2)} Let $A \to B$ be any map such that $B$ is reduced and
flat over $A$, and let $R_B = B \otimes_A R$. Then $\cJ(R_B/B)$
is defined,
and $\cJ(R_B/B) = \cJ(R/A)R_B$.  Note that there
is an induced map of total quotient rings $T \to T'$, and
it follows easily that $T' \otimes R_B$ is geometrically reduced and
equidimensional of dimension $d$. Also note
that $R_B$ is torsion-free over $B$:  $R$ is a directed union
of finitely generated $A$-submodules that are embeddable in
free $A$-modules, and since $B$ is flat over $A$ this property
is preserved by $B \otimes_A \blank$.

\item{(3)} In particular, we may apply (2) whenever $B$ is any
localization of $A$, or if $A$ is a field and $B$ is any
extension field.

\item{(4)} If $\cJ(R/A)$ is defined and $S$ is the localization
of $R$ at one element $f$, i.e., $S = R_f$,  then
$\cJ(S/A)$ is defined and equal  to $\cJ(R/A)S_f$.  Note
that if we have a presentation of $R$ such that $g(\vect x n)$ maps
to $f$,  then we get a presentation of $S$ by using one additional
variable $y$ and one additional generator for the denominator
ideal, $yg-1$,  and the calculation is then routine.

\item{(5)} Given ring extensions $A \to R$ and $R \to S$ such that
$\cJ(R/A)$ and $\cJ(S/R)$ are defined with $\cT(A) \otimes_A R$
of dimension $d$ and $\cT(R)\otimes_R S$ of dimension $d'$,  then
$\cJ(S/A)$ is defined,  $\cT(A)\otimes_AS$ has dimension
$d+d'$,  and $\cJ(S/A) \supseteq \cJ(S/R)\cJ(R/A)$. (Certainly,  $S$
is torsion-free over $A$.  The statements about being
reduced or geometrically reduced and about dimension can be checked
after tensoring with $T$, and we may, in fact, assume that
$A = T$ is a field.  The verifications are now straightforward.  For
example, the statement about dimension can be verified, for
each component of $S$, using the additivity of
transcendence degree.    For the statement about products of Jacobians
we can take a presentation of $S$ over $A$ of the form
$$
A[\vect x n, \vect y t]/( f_1(\ux),\,\ldots,\,f_m(\ux),
g_1(\ux,\uy), \,\ldots,\,g_s(\ux,\uy))
$$
where $\ux = \vect x n$ and $\uy = \vect y t$.  Here, we can
assume that
$$
A[\vect x n]/( f_1(\ux),\,\ldots,\,f_m(\ux) )
$$
is a presentation of $R$ over $A$ (let $M$ be the Jacobian
matrix with entries mapped to $R$), and that the images of
the $g_k$ in $R[\vect y t]$ may be used to give a presentation
of $S$ over $R$ (let $N$ be the corresponding Jacobian matrix  with
entries mapped to $S$).     Then the Jacobian matrix for $S$ over $A$ for
this presentation
with entries mapped to $S$  has the block form
$\pmatrix M & U\\ 0 & N\endpmatrix$.  Given  $n-d$ rows and columns
of $M$ (corresponding to the choice of a minor) and
$t-d'$ rows and columns of $N$, we get $(n+t)-(d+d')$ rows and
columns of this block matrix, and the determinant of the minor
they determine is the product of
the determinants of minors chosen from $M$ and $N$.)
\par
}
\enddemo

\demb{Discussion 3.2: test elements}
An element $c \in R\0$ is called {\it a test element} if, whenever $M$ is
a finitely generated $R$-module and $N \inc M$ is a submodule, then
$u \in M$ is in the tight closure of $N$ if and only if for
all $q = p^e$,  $cu^q \in N\q$ (the image of $\F^e(N) \to \F^e(M))$.
Thus, if the ring has a test element, it ``works" in any tight closure
test where some choice of $c \in R\0$ ``works."  Test elements
are also characterized as the elements of $R\0$ that annihilate
$N^*/N$ for all submodules $N$ of all finitely generated modules $M$.

A test element is called {\it locally stable} if its image in every local
ring of $R$ is a test element (this implies that it is a locally stable
test element in every localization of $R$ at any multiplicative
system).  A test element is called {\it completely stable} if
its image in the completion of each local ring of $R$ is
a test element:  a completely stable test element is easily
seen to be locally stable.  We refer the reader to
[HH1],
[HH2, \S6 and \S8], and [HH3, \S6] for more information about
test elements,  For the moment we simply want to note that
if $R$ is any reduced ring essentially of
finite type over an excellent local ring, then $R$ has a test element.
In fact, if  $c$ is any element of $R\0$ such that $R_c$ is regular
(and such elements always exist if $R$ is excellent and reduced),
then  $c$ has a power that is a completely stable test element.
This follows from Theorem (6.1a) of [HH3]. \enddemo

\demb{3.3 Discussion} When $R$ is a reduced ring of positive
prime characteristic $p$ and
$q = p^e$, we
write $R\qq$ for the unique reduced $R$-algebra obtained by
adjoining $q\,$th roots for all elements of $R$.  Thus, there is
a commutative diagram:
$$\CD  R @> \iota>> R^{1/q} \\
      @V 1_R VV   @VV \phi V \\
       R @>> F^e > R \endCD
       $$
where $\iota$ is an inclusion map and
$\phi(s) = s^q$.  We write $R\8$ for the
increasing union of the rings $R^{1/q}$.  The following result
is a variant of the results of \S1.5 of [HH6]:  the
differences from what is done in [HH6] are discussed in the proof. \enddemo

\proclaim{Theorem 3.4} Let $R$ be a geometrically reduced equidimensional
affine algebra of dimension $d$ over a field $K$ of positive
prime characteristic $p$. Let $t$ be an indeterminate over
$K$, let $L = K(t)$,  and let $R_L = L \otimes_K R$.
Let $J' = J_{R_L/L}$ be the
Jacobian ideal of $R_L$ over $L$,  which is evidently $\cJ(R/K)R_L$.
There are always elements of $\cJ(R/K)$ in $R\0$ (so that
$\cJ(R/K)$ is generated by such elements), and these are
completely stable test elements for $R$.  Moreover, $J'$
is generated by elements $c$ such that
{\pa
\item{\hbox{\rm ($*$)}} There is a regular subring $A$ of $R$
(depending on $c$),
in fact, a polynomial ring over $L$, such that $R_L$ is
module-finite and generically \'etale over $A$ and such that
for every $q = p^e$,   $cR\qq \inc A\qq[R]$:  moreover,
$A\qq[R] \cong A\qq\otimes_A R$ is $R$-flat for every $q$.\par
}
\endproclaim
\demb{Proof} We note that, in essence, all of this is established in
the proof of (1.5.5) of [HH6].  The fact that $\cJ(R/K)$ is not
contained in a minimal prime of $R$ follows from the fact
that $R$  is geometrically reduced.  The statement about completely
stable  test
elements is proved in (1.5.5) of [HH6] (there is an unnecessary
additional hypothesis in [HH6] that  $R$ be a domain --- we
discuss below why this can be removed).

The infinite field $L$ is needed so as
to be able to map a polynomial ring, say in $n$ variables,
onto $R$ in such a way
that the variables, after a suitable linear change of coordinates,
are in sufficiently general position.  Then, $R$ will be a module-finite
generically \'etale extension of any polynomial subring $A$ generated by $d$
of these variables, and it follows that every size $n-d$ minor
$c$ of the matrix
occurs in a Jacobian ideal $\cJ(R/A)$ (the notation agrees with
that used in \S(1.5.2) of [HH6]),
and so multiplies $R\qq$ into $A\qq[R]$ as a corollary of
the Lipman-Sathaye Jacobian theorem [LS].  There is one point that
needs a comment:  the Lipman-Sathaye theorem as given in [LS]
assumes that the ring $R$ is a domain, and because of this the result
in [HH6] is also stated with a domain hypothesis for $R$.  However, the
Lipman-Sathaye theorem
is valid in the reduced equidimensional case:  the needed result
is given in [Ho].  Finally,
we note that the isomorphism $A\qq[R] \cong A\qq\otimes_AR$
is proved in [HH2], Lemma (6.4), and since $A\qq$ is flat over $A$ (because
$A$ is regular: cf.\ [Kunz]) the
result follows. \qed\enddemo

\proclaim{Lemma 3.5} Let $R$ be a reduced Noetherian ring of
positive prime characteristic $p$,  let $c \in R$, and suppose that for
every power $q$ of $p$ there is an $R$-flat submodule $N_q$ of
$R\qq$  such that $cR\qq \inc N_q$.  Let $W$ be a multiplicative
system in $R$ and let $I$ be an ideal of $R$ that is contracted
with respect to $W$.  Let $I_q$ denote the contraction of
$I\q R_W$ to $R$.
Then for every $q = p^e$,  $c^qI_q \inc I\q$.
\endproclaim
\demb{Proof}  Since $N_q$ is $R$-flat,   $(R/I) \otimes _R N_q$
is $(R/I)$-flat.   Since  the elements of $W$ are nonzerodivisors
in $R/I$,  it follows that they are not zerodivisors on
$N_q/IN_q$.  If $u \in I_q$ we can choose $f \in W$ such that
$fu \in I\q$ and then $f^qu \in I\q$  as well.  Taking $q\,$th roots
we find that  $fu\qq \in IR\qq$,  and multiplying by $c$ gives
that $cfu\qq \in I(cR\qq) \inc IN_q$.  But $f$ is not a
zerodivisor on $N_q/IN_q$ (note that $cu\qq \in N_q$) and so
$cu\qq \in IN_q \inc IR\qq$.  Taking $q\,$th powers yields
that $c^qu \in I\q$,  as required.  \qed\enddemo

\proclaim{Lemma 3.6}  Let $R$ be a geometrically reduced equidimensional
$K$-algebra finitely generated over a field $K$ of
positive prime characteristic $p$.  Let $I$ be any ideal of $R$,
and let $W$ be a multiplicative system consisting of nonzerodivisors
modulo $I$.  Let $I_q$ be the contraction of $I\q R_W$ to $R$.
Then for every $q  = p^e$,  $\cJ(R/K)\q I_q \inc I\q$.
\endproclaim
\demb{Proof}  Let $L$ be as in Theorem 3.4.  After a flat base change
from $R$ to $R_L$ the image of $W$ still consists of nonzerodivisors
on $IR_L$,  and since $J_{R_L/L} = \cJ(R/K)R_L$ is generated
by elements $c$ satisfying the condition on $c$ in the hypothesis
of Lemma 3.5, if
$I_q'$ denotes the contraction of $I\q(R_L)_W$ to $R_L$, we have
that $J_{R_L/L}I_q' \inc I\q R_L$ and so $\cJ(R/K)I_q' \inc I\q R_L$.
Since $I_q \inc I_q'$,  it follows that $\cJ(R/K)I_q \inc
(I\q R_L)\cap R = I\q$, since $R_L$ is faithfully flat over $R$. \qed
\enddemo

\proclaim{Theorem 3.7} Let $R$ be a geometrically reduced equidimensional
$K$-algebra finitely generated over a field $K$ of positive prime
characteristic $p$.
Let $I$ be any ideal such that for every prime ideal $Q$ of
$R$, $IR_Q$ either contains a nonzerodivisor or else is (0) (i.e.,
$I$ is locally generically free).
Let $h$ be the largest analytic spread of $IR_P$ as $P$ runs through the
associated primes of
$I$.  Let $J = \cJ(R/K)$ be the Jacobian ideal.  Then for every
positive integer $n$ we have that
{\pa
\part{a} $J^n I\sy hn  \inc (I^n)^*$ (tight closure).
\part{b} $J^{n+1} I\sy hn \inc I^n$. \par}

More generally, for every nonnegative integer $k$ and positive integer $n$
we have that

{\pa
\part{a$'$} $J^n I\sy hn+kn \inc \bigl( (I\sy k+1 )^n \bigr)^*$
(tight closure).
\part{b$'$} $J^{n+1} I\sy hn+kn \inc (I\sy k+1 )^n$. \par
}
\endproclaim
\demb{Proof}  We have stated parts (a) and (b) separately for emphasis,
but evidently it suffices to prove the more general statements
(a$'$) and (b$'$).  Since $J$ consists of test elements it multiplies
the tight closure of any ideal into the ideal.  Thus, (b$'$) follows
from (a$'$) by multiplying by $J$, and it will suffice to prove (a$'$).

It suffices to prove the result for each connected component
of Spec $R$:  tight closures may be computed componentwise, and
passing to the component can be achieved by localizing at an idempotent ---
since formation of the Jacobian ideal commutes with localization, the new
Jacobian ideal is just the expansion of the original to the factor ring
corresponding to the component.

Thus, we may assume, by Lemma 2.5, that $I$ is either $(0)$ or else
contains a nonzerodivisor.  In the case where $I = (0)$ there is
nothing to prove.  If $I$ contains a nonzerodivisor then this is also
true for all powers of $I$.
Let $u \in I\sy hn+kn $.  Let $s$ be as in Lemma 2.4(b).

We must show that
$J^nu \inc \bigl( (I\sy k+1 )^n \bigr)^*$. Let $W$ be the
complement of the union of the associated primes of $I$.
For any $q = p^e$ we may
write $q = an + r$
with $a$ a nonnegative integer and $0 \leq r < n$,  and then
by Lemma 2.4(b) we have that
$$
I^{s + (h+k)(n-1)}u^a \inc \bigl((I\sy k+1 )\q\bigr)^W
$$
where the superscript $^W$ indicates expansion to $R_W$ followed
by contraction to $R$.
Since $I\sy k+1 $ is contracted
with respect to $R_W$, we may use Lemma 3.6 to conclude that
$$
J\q \bigl((I\sy k+1 )\q\bigr)^W \inc (I\sy k+1 )\q ,
$$
and so we have  that
$$
J\q  I^{s + (h+k)(n-1)}u^a \inc (I\sy k+1 )\q .
$$
Taking $n\,$th powers and abbreviating  $b = sn + (h+k)(n-1)n$ we have
that
$$
I^b (J^n)\q u^{an} \inc \bigl((I\sy k+1 )^n\bigr)\q
$$
for all $q$ and since $q \geq an$ this yields
$$
I^b (J^nu)\q \inc \bigl((I\sy k+1 )^n\bigr)\q
$$
for all $q$.  Let $d$ be a fixed nonzerodivisor in $I^b$ (note that
$b$ does not depend on $q$).  The condition that
$$
d(J^nu)\q \inc \bigl((I\sy k+1 )^n\bigr)\q
$$
tells us precisely that $J^n u \inc \bigl((I\sy k+1 )^n\bigr)^*$,
as required. \qed \enddemo

\demb{Example 3.8} Consider the ring $R = K[x,y,z]/(xy-z^n)$.
The Jacobian ideal $J= \cJ(R/K)$, if
$n$ is a unit, is $(x,y,z^{n-1})R$. Let $P = (y,z)$. Then
$h = 1$, and  $y\in P\sy n $.
Then $J^{n-1}$ multiplies $y$ into $P^n$ but no smaller
power does since $y\in J$ and
$y^{n-2}y\notin P^n$. This suggests that the result in 3.7 (take
$I = P$) is close
to best possible.  We do not know, however, whether the exponent
$n$ used in parts (a) and (a$'$) can be replaced by $n-1$ in general.
Of course, if so, then the exponent $n+1$ can be replaced by
$n$ in parts (b) and (b$'$).  \enddemo
\bigskip
\bigskip
{\chbf
\centerline{\hdbf 4. The equal characteristic zero case} \smallskip
}
\bigskip
We now give the extensions of the various positive characteristic
results to the equal characteristic case.  As mentioned in the
introduction, the notion of tight
closure that we use here is that of equational tight closure
from [HH6, \S\S 3.4.3--4].  The main results of this section
are contained in Theorem 4.4 below.  We need to do some groundwork
before we can prove that theorem, however.
The proof of the main results depends on three steps:  one is to
localize and complete, the second is to descend from the complete
case to the affine case, and the third is to use reduction to
positive characteristic in the affine case.  The second step
is based on the following result from [AR]:

\proclaim{Theorem 4.1} Let $K$ denote either a field or
an excellent discrete valuation ring.
Let $T=K[[x_1,\ldots,x_n]]$ be the formal power series ring in $n$ variables
over $K$.
Then every $K$-algebra homomorphism of a finitely generated $K$-algebra $R$
to $T$ factors $R \to S \to T$ where the maps are $K$-algebra homomorphisms
and $S$ has the form $(K[x_1,\ldots,x_n, y_1,\ldots,y_t]_m)\2$, where the
$x_i$ are as above, the $x_i$ and $y_j$ are algebraically independent
elements, over
$K$, of the maximal ideal of $T$, $m$ is the ideal of the polynomial ring
$K[x,y]$ generated by $(x,y)$ and, if $K$ is a DVR, by the generator of the
maximal ideal of $K$, and $\,\2$ denotes Henselization. \qed
\endproclaim

This is a special case of general N\'eron desingularization
(cf.\ [Po1], [Po2], [Og], [Swan]), but the argument is simpler
in this case (we note that [Swan] has removed any possible
doubt about the validity of the general theorem --- however, we
only need the result of [AR]).   In [HH6] this is used to prove
the following result, which is Theorem (3.5.1) there:

\proclaim{Theorem 4.2} Let $K$ be a field of characteristic zero and let
$(S,m,L)$ be a complete local ring that is a $K$-algebra.
Assume that $S$ is equidimensional and unmixed.

Suppose that $R_0$ is a subring of $S$ that is finitely generated as a
$K$-algebra.
We also assume given finitely many sequences of elements $\{z_t^{(i)}\}{}_t$
in $R_0$, each of which is part of a system of parameters for $S$.

Then there is a finitely generated $K$-algebra $R$ such that the
homomorphism  $R_0 \inj S$  factors $R_0 \inj R \to S$ and such that
the following conditions are satisfied:\smallskip

{\pa
\part{1} $R$ is biequidimensional.

\part{2} The image of each sequence $\{z_t^{(i)}\}_t$ in $R$ is a sequence of
{\it strong} parameters:  this means that after localization and completion
at any prime that contains them, they form part of a system of
parameters modulo every minimal prime.

\part{3} If $\m$ is the contraction of $m$ to $R$, then
$\dim R_\m - \depth R_\m = \dim S - \depth S$.
In particular, $R_\m$ is Cohen-Macaulay iff $S$ is Cohen-Macaulay.

\part{4} If $S$ is a reduced (respectively, a domain) then so is $R$.\par
}
\smallskip

(N.B. In general, $\dim R_m$ is substantially bigger than $\dim S$.) \qed
\endproclaim

In the sequel we need a version of this result in which the
equidimensionality of the ring is not assumed.  Moreover, we need
to keep track of some complexes of modules, bounds on analytic
spreads after localization at associated primes,
the fact that one ideal is a certain symbolic
power of another, and so forth.   The following result suffices:

\proclaim{Theorem 4.3} Let $K$ be a field of characteristic zero and let
$(S,m,L)$ be a complete local ring that is a $K$-algebra.
Assume given each of the following:
{\pa
\part{1} Finitely many finitely generated $S$-modules
with specific finite presentations,
finitely many maps of these $S$-modules with specific
presentations of the maps, and  finitely many
specified equalities among the compositions of these maps.

\part{2}
Using  the modules and maps in (1), finitely many short exact sequences.
Finitely many finite complexes with specified homology. Finitely
many instances in which one of the specified modules is
identified with a submodule of another.  Finitely many
instances in which one of the modules is specified to
be the intersection of finitely many of the others, where
all are submodules of a given specified module.

\part{3} Finitely many ideals of $S$ with specified
generators.

\part{4} From among the ideals in (3), a finite subset
with a finite set of associated primes of specified heights,
and a finite subset such that the maximum analytic spread after localizing
at an associated prime has a given bound.

\part{5} From among the ideals in (4), finitely many choices of
$I$, $I'$, such that $I$ and $I'$ have the same
associated primes.  Also, finitely
many choices of $I$ and $I'$ and an integer $k$ such that $I' = I\sy k $.

\part{6} From among the ideals and modules given in (1) and (3),
finitely many pairs $M$, $I$ such that $IM = 0$.

\part{7} A finite set of finite sequences in $S$, each of which
is specified to be part of a system of parameters for $S$,
and a  finite set of sequences each of which is specified to be a
regular sequence on a given one of the given modules.

\part{8}
A finitely generated $K$-subalgebra $R_0 \inc S$  so large
that it contains all the entries needed for the
presentations of the modules and maps in (1),
all of the specified
generators of the ideals in (3), and the elements of the
sequences in (7),  so that we may
view all of the given modules, maps, sequences, and ideals as arising
from corresponding ones over $R_0$  either by tensoring, taking images, or
expanding ideals.  \par
}  \smallskip
Then there is a finitely generated $K$-algebra $R$ such that the
homomorphism  $R_0 \inj S$  factors $R_0 \inj R \to S$ and such that
the following conditions are satisfied:\smallskip

{\pa
\part{a} The specified presentations of maps of modules are maps of
modules over $R$, and the specified exacts sequences of modules, descended
to $R$ by tensoring up from the their counterparts over $R_0$, are
exact.   All of the other specified relations among the given
modules and ideals continue to hold after descent, including
specifications of the homology of a given complex and specifications
that a certain submodule (or ideal) be a finite intersection of
finitely many given other submodules (or ideals).   Likewise,
the specification that a certain ideal be the annihilator of a certain
module can be preserved.

\part{b} The image (under the map $R_0 \to R$) of each set of
elements that is part of a
system of parameters for $S$ has height equal to its length.
The image of each regular sequence
on a specified module is a regular sequence on the corresponding
module over $R$.

\part{c} The specified ideals, descended to $R$ by expanding
their counterparts over $R_0$, are unmixed when the original
ideals are.  For a specified ideal $I$, the greatest number
of generators and the greatest analytic spread after localization at an
associated prime do not
increase.    Moreover, for the given choices of $I$, $I'$, $k$
such that $I' = I\sy k $,  this remains true after descent to $R$.

\part{d} $R$ is regular if $S$ is.
\par
}
\smallskip

(N.B. In general, $\dim R_m$ is substantially bigger than $\dim S$.)
\endproclaim

\demb{Proof} If $S$ is regular we apply the Artin-Rotthaus theorem
(4.1) directly to the power series ring $S = L[[\vect x n]]$, where
the coefficient field $L$ has been chosen to contain $K$.  We first
solve the problem over $L$ and then  descend to $K$.  As the latter
step is routine, we shall simply treat the case $L = K$.

We are free to enlarge
$R_0$ repeatedly, and so may assume that $\vect x n \in R_0$.
Since the Henselization of a local ring
is a direct limit of finitely presented \'etale extensions,
we have that $S$ is a filtered inductive limit of  regular rings
$R$ of finite type over $L$ with maps
$R_0 \to R \to S$ such that  $R$ is smooth over
$L[\vect x n]$ and such that the $x_i$ form a permutable regular
sequence in $R$.  In this case, we can keep track of whether
a sequence of elements is part of a system of parameters by
extending it to a full system of parameters, say, $\vect y n$.
There will be equations expressing a power of every $y_j$ as
a linear combination of the $x_i$ and conversely.  We may
enlarge $R_0$ so that all these equations hold in $R_0$.
Conditions such as having a specified ideal kill a specified
module are likewise expressible equationally and can
be guaranteed by enlarging $R_0$.  We shall leave many
straightforwarded details to the reader, noting that
the proof of (3.5.1) in [HH6] is given in great detail and
is a very similar kind of argument.  We focus here only
on some critical issues.

One can keep track of short exact sequences
by using an exact sequence of finite free resolutions over $R$.
The rows will be split exact.  To guarantee that a finite
free resolution stays a resolution one keeps track of all
the matrices.  $R_0$ is enlarged to contain all their entries.
The condition that one has a complex is equational, and so
is the condition that the determinantal ranks be preserved.
By the result of [BE],  one only needs to guarantee that
the largest nonvanishing ideals of minors have specified
depths, i.e., that each contains a subset, of a certain
specified size,  of a system of
parameters for the ring, and we may apply the discussion
of the preceding paragraph.

This enables one also to keep
track of finite complexes with specified homology (express all
the conditions by using suitable
short exact sequences) and of finite
intersections as well: e.g., the intersection of $N_1$ and
$N_2$ within $N$ may be characterized as the kernel
of the map $N \to N/N_1 \oplus N/N_2$.  The annihilator $I$
of a single element $u$ of a module $N$ may be characterized
by an injection $R/I \to N$ carrying the image of $1$ to $u$,
and the annihilator of $N$ may be characterized as the intersection
of the annihilators of specified generators of $N$.

One can preserve depths of modules
and, hence, regular sequences by expressing them in terms
of the vanishing of  Koszul homology.

One can keep the associated primes of an ideal having
specified heights as follows:  keep track of
its entire primary decomposition, preserving the fact that
components intersect to give the ideal.  To preserve the
relation between an ideal $I$ primary to $P$ and and
$P$,  note that $S/I$ has a filtration by modules that
are embeddable in finitely generated free $(S/P)$-modules,
and this can be preserved.  Now, the height property
can be preserved by keeping $P$ height unmixed in the descent:
we do not need to keep $P$ prime.  This can be achieved
by writing $P$ in the form
$$
(\vect f d)R\col R gR
$$
where the
elements $f_i$ are part of a system of parameters for $S$.
Note that the number of primary components may increase,
but the largest height of an associated prime does not.
The same idea can be used to ensure that two specified ideals
$I$, $I'$ that have the same associated primes continue to
do so.

When $S$ is not necessarily regular write it as $T/J_0$
where $T$ is regular, and transfer the problem to $T$ (while
keeping track of $J_0$).  Ideals of $S$ correspond to
ideals of $T$ that contain $J_0$, and $R$-modules to $T$-modules
that are killed by $J_0$.  For example, to maintain a specific symbolic
power relationship, one may suppose that one of the ideals
is $I \supseteq J_0$ and that it has a certain set $\cS$ of associated
primes while the other has the form
$(I^k + J_0) \Col g$ with $g$ a nonzerodivisor on $I$ and that it has
as its associated primes a certain subset of $\cS$.
We have seen that all this can be preserved
while descending.

Finally, we want to explain how to preserve the condition that
the maximum analytic spread of $I$ after localizing at an
associated prime of $I$ be at most $h$.  Again, we think
of $S$ as $T/J_0$ where $T$ is regular.  Call the associated
primes $\vect P s$.  Then for each $P_i$,  $IS_{P_i}$ is integral
over $IR_{P_i}$ after localizing at $P_i$.  Thus,  for each $P_i$
we can choose an element $v = v_i$  not a zerodivisor on $P_i$
such that $IR_v$ is integral over $I_0R_v$,  where $I_0$ is
generated by at most $h$ elements of $R$.  (In the equal characteristic
0 case, the residue field is always infinite.)  After clearing
denominators by multiplying by a power of $v$,  for each generator $r$
of $I$  we get an equation
$$
v^Nr^t + i_1r^{t-1} +\,\cdots\,+i_kr^k+\,\cdots\,+i_{t-1}r + i_t = j ,
$$
with $i_k \in I_0^k$ and with $j \in J_0$.
We can preserve all this by placing all of the
needed elements in $R_0$.  As we descend, the $P_i$ are replaced
by unmixed ideals, while each $v_i$ is kept a
nonzerodivisor on the descended version of $P_i$.
By including sufficiently many coefficients in $R_0$ we can
preserve that every $i_k \in I_0^k$ where $I_0$ is generated
by at most $h$ elements.
 Any associated prime
of the descended version of $I$ will be an associated prime of one of
the descended $P_i$, and so will fail to contain the corresponding
$v_i$.  Thus, after descent, when one localizes at an associated
prime of the descended version of  $I$, at least one of the $v_i$
becomes invertible, and
it follows that the descended version of $I$ becomes integrally
dependent on an ideal
with at most $h$ generators.
\qed\enddemo

\proclaim{Theorem 4.4} Let $R$ be Noetherian
ring containing a field of characteristic $0$.  Let $I$ be any
ideal of $R$, and let $h$ be the largest analytic spread
of $IR_P$ for $P$ an associated prime
of $I$.
{\pa
\part{a} If $R$ is regular, $I\sy hn \inc I^n$ for all positive integers $n$.
More generally,
$I\sy hn+kn \inc (I\sy k+1 )^n$
for all positive integers $n$ and  nonnegative integers $k$.

\part{b} If $I$ has finite projective
dimension then $I\sy hn  \inc (I^n)^*$ for all positive integers $n$.

\part{c} If $R$ is affine and equidimensional over
a field $K$, and locally $I$ is either $0$ or contains a nonzerodivisor,
then with $J = \cJ(R/K)$,  for every nonnegative integer $k$ and
positive integer $n$ we have $J^nI\sy hn+kn \inc ((I\sy k+1 )^n)^*$ and
$J^{n+1}I\sy hn+kn \inc (I\sy k+1 )^n$.   In particular,
$J^nI\sy hn \inc (I^n)^*$ and $J^{n+1}I\sy hn  \inc I^n$ for all $n$.
\par
 }
\endproclaim
\demb{Proof} We first prove (c), and at the same time we prove
(b) for finitely generated algebras over a field $K$.  We use
the standard descent theory of Chapter 2 of [HH6] to replace
the field $K$ by a finitely generated $Z$-subalgebra $A$, so
that we have a counterexample in an affine algebra $R_A$ over
$A$ with $R_A \inc R$ and $R \cong K \otimes_A R_A$.
In particular, $R_A$ will be reduced.
In doing so we descend $I$ to an ideal $I_A$ of $R_A$
as well as the ideals and their prime radicals in its primary
decomposition.  We have an element $u_A$ that fails  to satisfy
the containment we are trying to prove.  In the regular case,
we can localize at a nonzero element of $A$ to make $R_A$ smooth
over $A$.  In either case, we can localize at a nonzero element
of $A$ to make $A$ smooth over $\Z$.  Since $\cJ(R_A/A)\cJ(A/\Z)
\inc \cJ(R_A/\Z)$ and since $\cJ(A/\Z) = A$ when $A$ is smooth
over $\Z$,  we see that we may assume that $\cJ(R_A/A) \inc
\cJ(R_A/\Z)$,  which means that we can work over $\Z$ instead of
$A$.  The result now follows from the fact that, for almost all
fibers, the containment
holds for the map $\Z \to R_A$ after passing to fibers over closed
points of $\Spec \Z$.  Notice that as we pass to fibers
$\kappa \to R_\kappa$
we may assume that each minimal prime $P_A$ of $I_A$ becomes
a radical ideal whose
minimal primes in $R_{\kappa}$ are all of the same height as
the original.  Thus, in the fiber, the primary decomposition
of $I_{\kappa}$ may have more components, but each of these
will be obtained from the image of one of the original components
by localization.  The biggest analytic spread after localizing at
an associated prime will not change.
It follows in both parts that we have the required containment in a
tight closure.  In the regular case, we have that all ideals are
tightly closed.

We now consider the general case for (a) and (b).  The problem
in each part reduces to the local case:  note that it suffices
to check whether an element is in a tight closure locally after
completion.

One may then complete:  although $I\Rc$ may have more associated
primes, Discussion 2.3(c) shows that the biggest analytic
spread as one localizes at these  cannot increase.
Note that once we have $\Ic\sy hn \inc (\Ic^n)^*$  in $\Rc$,
it follows that $I \sy hn \inc \bigl((I^n)\Rc\bigr)^*$  since
$\Rc$ is flat over $R$,  and this implies that $I \sy hn \inc (I^n)^*$,
which is $I^n$  when $R$ is regular.

In the regular case,
next note that $\Ic\sy k+1 = I\sy k+1 \Rc$, so that
$(\Ic\sy k+1 )^n = (I\sy k+1 )^n\Rc$. (The associated
primes of  $\Rc/I\sy k+1 \Rc$ are among those associated
to $\Rc/P\Rc$ for some associated prime $P$ of $I$, by Proposition 15
in IV B.4. of [Se], since any associated prime of $I\sy k+1 $
must be an associated prime of $I$, and by another application of
Proposition 15
in IV B.4. of [Se] these
in turn are associated primes of $I\Rc$.)  Thus, we get
$$
I\sy hn+kn \inc \Ic\sy hn+kn \inc (\Ic\sy k+1 )^n =  (I\sy k+1 )^n\Rc
$$
and so
$$
I\sy hn+kn \inc (I\sy k+1 )^n\Rc \cap R = (I\sy k+1 )^n ,
$$
as required, by the faithful flatness of $\Rc$ over $R$.

Using Theorem 4.3 above
one may then descend to a suitable affine algebra over a coefficient
field for the complete local ring, and the results follow from
what we have already proved in the affine case. \qed\enddemo
\bigskip
\bigskip
{\chbf
\centerline{\hdbf 5. Questions}
}
\bigskip
Evidently, if we fix an ideal $I$ in a Noetherian ring
$R$, for every integer $N$ there is a least integer $g(N) \in \N$
such that $I \sy g(N) \inc I^N$, and there are clearly deep results
about the behavior of $g(N)/N$.  Our result that in equicharacteristic
regular rings (or when $I$ has finite projective dimension),
$g(N)/N$ is bounded by
the largest height of an associated prime of $I$
(or the largest analytic spread of $IR_P$ for $P$ an
associated prime of $I$)  might be
improved in any number of ways. \medskip

One possibility is to study the case where the ring
need not be  regular (and $I$ does not have
finite projective dimension).  Note that our results here for this
case, involving the Jacobian ideal, do not directly give information
about  the question raised just above.  (We  mention again that
we do not know whether the exponents used on the Jacobian ideal
in Theorem 3.7 and Theorem 4.4(c) are best possible:  it may
be possible to decrease the exponent by 1.)      \medskip

We do not know what the situation is
in mixed characteristic regular rings.   But even in
equicharacteristic regular rings there may be better bounds that
make use of additional information about $I$.  Notice, for example,
that the height is never the best bound when $I$ is $m$-primary
in a regular local ring $(R,\,m)$,  since then $I\sy n  = I^n$
for all $n$. \medskip

It is not clear what the best bound is in equicharacteristic
regular rings even for primes of codimension 2.  \medskip

In quite a different direction, we observe that
there have been several instances in which the theory of
multiplier ideals and tight closure theory have either interacted,
or have been used to prove similar results. E.g., tight closure
can be used to prove the Brian\c con-Skoda theorem (cf.\  [HH2, \S5]),
as can the theory of multiplier ideals. (Cf.\ [L2], where these
are called {\it adjoint} ideals. This is also done
implicitly in [EL] and explicitly in [Laz, \S10]).  The connection
between the multiplier ideal and the test ideal of tight closure theory
is explored in [Sm2]. It would be desirable to understand fully
the underlying reasons for this connection.

} 

{
\smc
\baselineskip = 10 pt
\settabs 7 \columns
\quad\bigskip
\+Department of Mathematics           &&&&Department of Mathematics\cr
\+University of Michigan              &&&&University of Kansas\cr
\+Ann Arbor, MI 48109--1109           &&&&Lawrence, KS 66045\cr
\+USA                                 &&&&USA\cr
\smallskip
\+E-mail:                             &&&&E-mail:\cr
\vskip 1.5 pt plus .5 pt minus .5 pt
{\rm
\+hochster\@umich.edu                 &&&&huneke\@math.ukans.edu\cr
}
}

\enddocument